\documentclass[12pt,reqno]{amsart}
\usepackage{amsthm}
\usepackage{cases}
\usepackage{amscd}
\usepackage{amsfonts}
\usepackage[dvipdfmx]{graphicx}
\usepackage{makecell}
\usepackage{amssymb}
\usepackage{amsmath}
\usepackage{booktabs}
\usepackage{amsthm}
\newtheorem{example}{Example}

\usepackage{subfigure}
\usepackage{algorithm}
\usepackage{algpseudocode}
\renewcommand{\algorithmicrequire}

\usepackage{cleveref}
\usepackage{color}
\usepackage{datetime}

\allowdisplaybreaks

\makeatletter
\newcommand\figcaption{\def\@captype{figure}\caption}
\newcommand\tabcaption{\def\@captype{table}\caption}
\makeatother

\usepackage{geometry}
\usepackage{algorithm}
\usepackage{multirow}
\usepackage{mathrsfs}
\usepackage{ulem}

\newcommand{\beq}{\begin{equation}}
\newcommand{\eeq}{\end{equation}}
\newcommand{\bea}{\begin{eqnarray}}
\newcommand{\eea}{\end{eqnarray}}
\newcommand{\beas}{\begin{eqnarray*}}
\newcommand{\eeas}{\end{eqnarray*}}

\begin{document}
\title[DataInNet for solving PDEs]{Data-integrated neural networks for solving partial differential equations}

\author[J. Zheng, Y. Huang, N. Yi, Y. Yang]{Jiachun Zheng$^\dagger$, Yunqing Huang$^{\ddagger,*}$, Nianyu Yi$^{\S}$ and Yunlei Yang$^{\imath}$}

\address{$^\ddagger$ School of Mathematics and Computational Science, Xiangtan University, Xiangtan 411105, P.R.China }
\email{jczheng2022@126.com}

\address{$\ddagger$ National Center for Applied Mathematics in Hunan, Key Laboratory of Intelligent Computing \& Information Processing of Ministry of Education, Xiangtan University, Xiangtan 411105, Hunan, P.R.China} 
\email{huangyq@xtu.edu.cn}

\address{$^\S$ Hunan Key Laboratory for Computation and Simulation in Science and Engineering; School of Mathematics and Computational Science, Xiangtan University, Xiangtan 411105, P.R.China}
\email{yinianyu@xtu.edu.cn}

\address{$^\imath$ School of Mathematics and Statistics, Guizhou University, Guiyang, Guizhou 550025, P.R.China}
\email{ylyang5@gzu.edu.cn}



\begin{abstract}
In this work, we propose data-integrated neural networks (DataInNet) for solving partial differential equations (PDEs), offering a novel approach to leveraging data (e.g., source terms, initial conditions, and boundary conditions). The core of this work lies in the integration of data into a unified network framework. DataInNet comprises two subnetworks: a data integration neural network responsible for accommodating and fusing various types of data, and a fully connected neural network dedicated to learning the residual physical information not captured by the data integration neural network. This network architecture inherently excludes function classes that violate known physical constraints, thereby substantially narrowing the solution search space. Numerical experiments demonstrate that the proposed DataInNet delivers superior performance on challenging problems, such as the Helmholtz equation (relative \(L^2\) error: O(\(10^{-6}\))) and PDEs with high frequency solutions (relative \(L^2\) error: O(\(10^{-5}\))).
\end{abstract}


\thanks{$^*$ Corresponding author.}

\maketitle


\section{Introduction}\label{section-1}
In scientific computing, deep neural networks (DNNs) have been widely adopted for solving various partial differential equations (PDEs) \cite{1,2,3,4,5,6,7}. While notable progress has been made in PDE solving using DNNs \cite{8,9,10,11,12,13,14,15,16,17,18}, several critical challenges remain unaddressed, such as spectral bias (i.e., the frequency principle) \cite{19}. This bias is characterized by the tendency of DNNs to prioritize learning low frequency components over high frequency ones when fitting functions or solving PDEs. The existence of spectral bias has been experimentally validated in numerous studies \cite{20,21,22}, and Tancik et al. leveraged the neural tangent kernel \cite{23} to provide corresponding theoretical justification.
Inspired by insights into the neural tangent kernel, physics-informed neural networks (PINNs) incorporate multi-scale random Fourier feature mapping to enhance robustness and accuracy \cite{24,25}. Furthermore, PINNs leverage the derivatives of the solution to define the loss function-a design that may facilitate the convergence of high frequency components during training. Examples presented by Lu et al. \cite{21} confirm that differentiation promotes the convergence of high frequency components. 

Low numerical accuracy and training instability are the primary challenges for PINNs. To improve the approximation of PINNs, Lu et al. \cite{26} proposed a residual-based adaptive sampling strategy (RAR). RAR improves the approximation accuracy of PINNs by adding new training points with high residuals. In \cite{27}, Fan et al. developed discontinuity-removing physics-informed neural networks (DR-PINNs) for solving elliptic interface problems. DR-PINNs provide a high-precision approximation by using the geodesic acceleration Levenberg-Marquardt optimizer to minimize the loss function. Multi-stage neural networks \cite{28} divide the training process into different stages, where each stage employs a new network to fit the residuals from the previous stage. Numerical results show that multi-stage training for both regression tasks and PINNs yields approximations with prediction errors close to machine precision. Huang et al. \cite{29} proposed frequency-adaptive multi-scale DNNs (frequency-adaptive MscaleDNNs) for solving PDEs with high frequency solutions. The frequency-adaptive MscaleDNNs employ a multi-stage training strategy. In the first stage, MscaleDNNs are used to obtain a pre-trained solution. The frequency coefficient distribution of the pre-trained solution is obtained via Fourier transform, and these corresponding frequencies are then utilized to construct a new neural network with multiple subnetworks. This new network is trained to achieve a more accurate approximation, and this process is iterated until training is complete. Additionally, adaptive activation functions \cite{30}, adaptive weighting of loss functions \cite{31}, and learning rate annealing algorithm \cite{32} have also been proposed to improve the approximation of PINNs.

Although numerous PINN variants have been proposed for PDE solving, the data (e.g., source terms, initial conditions, and boundary conditions) are typically incorporated as constraints in the loss function. To address this limitation, we develop data-integrated neural networks (DataInNet) for PDE solution. The main contributions of this work are summarized as follows: 

\begin{itemize}
\item \textbf{Data integration architecture:} 
This architecture reformulates the conventional paradigm of PINNs, where source terms, initial conditions, and boundary conditions are merely formulated as soft constraints in the loss function. Instead, it innovatively integrates such physical data directly into the network architecture, providing a new reference for the utilization of physical prior information.

\item \textbf{Significant reduction of solution space search range:}  
By inherently excluding solutions that violate physical constraints, this structure enables the optimization process to focus more on feasible function classes. Consequently, it enhances solution accuracy, computational efficiency, and generalization ability, mitigating the optimization challenges commonly encountered in traditional PINNs.

\item \textbf{Dual-branch network structure:} 
Through a data integration network that structurally embeds known information (e.g., source terms, boundary conditions) and an auxiliary network responsible for learning remaining physical laws, our method achieves explicit separation and collaborative learning between known information and the information to be solved.
\end{itemize}

The rest of the paper is organized as follows: In section \ref{section-2}, we explain the network structure of data-integrated neural networks in detail. In section \ref{section-3}, the different PDEs are used to test the approximation ability of DataInNet. Finally, the advantages and inadequacies of DataInNet are summarized in section \ref{section-4}.

\section{DataInNet for solving PDEs}\label{section-2}
In this section, we first briefly introduce the PINNs method, then we 
propose the data-integrated neural networks for solving partial differential equations. 

\subsection{Physics-informed neural networks}
Considering the following problem:
\begin{equation}\label{2.1}
\left\{\begin{split}
&D u(\textbf{x},t)= f(\textbf{x},t),\quad \quad (\textbf{x},t)\in \Omega \times [0,T],\\
&u(\textbf{x},t)=B(\textbf{x},t),\,\;\,\,\quad \quad (\textbf{x},t)\in\partial\Omega \times [0,T],\\
&u(\textbf{x},0)=I(\textbf{x}),\,\,\,\;\quad\quad \quad \textbf{x}\in\Omega,\\
\end{split}\right.
\end{equation}
where $\Omega \in \mathbb{R}^{N}$ denotes the computational domain, and $D$ represents a differential operator containing temporal and spatial derivatives. The corresponding residual functions of problem (\ref{2.1}) are defined as:
\begin{equation}\label{2.2}
\begin{split}
&r(\textbf{x},t,\theta)=D u_{nn}(\textbf{x},t,\theta)- f(\textbf{x},t),\\
&r_{b}(\textbf{x}_{b},t_{b},\theta)=u_{nn}(\textbf{x}_{b},t_{b},\theta)-B(\textbf{x}_{b},t_{b}),\\
&r_{0}(\textbf{x}_{0},0,\theta)=u_{nn}(\textbf{x}_{0},0,\theta)-I(\textbf{x}_{0}),\\
\end{split}
\end{equation}
where $(\textbf{x},t)$, $(\textbf{x}_{b},t_{b})$ and $(\textbf{x}_{0},0)$ denote the interior point, the boundary point, and the initial point, respectively. $u_{nn}(\textbf{x},t,\theta)$ denotes a fully connected neural network with $M$ hidden layers defined as follows:
\begin{equation}\label{2.3}
\begin{aligned}
&\Phi^{0}(\bar{\textbf{x}})=\bar{\textbf{x}},\\
&\Phi^{m}(\bar{\textbf{x}})=\sigma(W^{m}\Phi^{m-1}(\bar{\textbf{x}}) + b^{m}),\qquad 1\leq m\leq M,\\
&u_{nn}(\bar{\textbf{x}},W)=W^{M+1}\Phi^{M}(\bar{\textbf{x}})+ b^{M+1},\\
\end{aligned}
\end{equation}
where $\bar{\textbf{x}}=\{\textbf{x},t\}$ and $\Phi^{m}(\bar{\textbf{x}})$ represent the input and hidden state of the neural network, respectively. $\sigma(\cdot)$ is a nonlinear activation function. $\theta=\{W^{m},b^{m}\}_{m=1}^{M+1}$ denote all the trainable parameters. Then, a neural network approximation for problem (\ref{2.1}) is obtained by minimizing the following loss function:
\begin{equation}\label{2.4}
\begin{split}
L(\textbf{x},t,\theta)&=\tau\frac{1}{N} \sum\limits_{i=1}^{N} r(\textbf{x}^{i},t^{i},\theta)^{2} + \tau_{b}\frac{1}{N_{b}}\sum\limits_{i=1}^{N_{b}} r_{b}(\textbf{x}_{b}^{i},t_{b}^{i},\theta)^{2}+ \tau_{0}\frac{1}{N_{0}}\sum\limits_{i=1}^{N_{0}} r_{0}(\textbf{x}_{0}^{i},0,\theta)^{2},\\
\end{split}
\end{equation}
where $\{\tau, \tau_{b}, \tau_{0}\}$ denote the weights of each of the components in the loss function.

\subsection{Data-integrated neural networks}
PINNs embed physical laws into the loss function as penalty terms and optimize the loss function to compel the neural network to approximately satisfy physical constraints. When dealing with various physical problems, neural networks often need to learn physical laws from scratch. This raises the question: if certain components of PDEs are known (e.g., source terms, initial conditions, and boundary conditions), why not directly embed such known components into the neural network architecture itself, instead of requiring the network to learn them from scratch? To address this issue, we develop data-integrated neural networks (DataInNet). 

DataInNet comprises a data integration neural network and an auxiliary neural network. The data integration neural network is specifically designed to accept and integrate various types of data, including source terms, initial conditions, and boundary conditions. The auxiliary neural network is a fully connected network tasked with learning the residual physical information not captured by the data integration neural network.
The forward propagation rule for data integration neural networks with $N^{*}$ modules are as follows: 
\begin{equation}\label{3.2}
\begin{aligned}
&H^{0}=\bar{\textbf{x}},\\
&H^{n}=\sigma(W^{n}A^{n}+b^{n}), \quad 1 \leq n \leq N^{*},\\
&u_{nn}^{d}(\bar{\textbf{x}})=W^{N^{*}+1}H^{N^{*}}+b^{N^{*}+1},\\
\end{aligned}
\end{equation}
where $A^{n}$ is defined as follows:
\begin{equation}\label{A}
\begin{aligned}
&G^{n}=W^{n}_{G}H^{n-1}+b_{G}^{n},\\
&Q^{n,i}=W^{n,i}_{Q}H^{n-1}+b_{Q}^{n,i},\\
&A^{n,i}=\sigma(G^{n}\odot(Q^{n,i}+\hbar^{i}(\bar{\textbf{x}}))),\quad i=f,B,I,\\
&A^{n} = \sum\alpha^{i}A^{n,i}.\\
\end{aligned}
\end{equation}
Here, $\odot$ denotes a element-wise multiplication, and $\sum\limits_{i=1}^{N}\alpha^{i}=1$. $\hbar^{i}(\bar{\textbf{x}})$ represents data:
\begin{equation}
\begin{split}
\hbar^{i}(\bar{\textbf{x}})=\left\{
	\begin{aligned}
	&f(\bar{\textbf{x}}),  \qquad\qquad\qquad i=f,\\
    &I(\bar{\textbf{x}}),  \qquad\qquad\qquad i=I,\\
    &\sum_{j=1}^{J} \alpha_{B}^{j} B^{j}(\bar{\textbf{x}}),  \,\,\qquad i=B,\\
	\end{aligned}
	\right.
\end{split}
\end{equation}
where $f(\bar{\textbf{x}})$, $I(\bar{\textbf{x}})$ represent the source term and initial condition, respectively. For PDEs defined on a polygon region, $B^{j}(\bar{\textbf{x}})$ denotes the boundary condition on the $j$-th edge, and $\alpha_{B}^{j}$ is the corresponding learnable weight. In each module $H^{n}$, we design a multi-branch parallel network architecture: each branch $Q^{n,i}$ is dedicated to incorporating specific data $\hbar^{i}(\bar{\textbf{x}})$, while the branch $G^{n}$ focuses on learning other intrinsic features of the input data $H^{n-1}$. Subsequently, the outputs of $G^{n}$ and $Q^{n,i}+\hbar^{i}(\bar{\textbf{x}})$ are fused via element-wise multiplication after dimension matching, enabling effective cross-feature interaction. Finally, a gating mechanism is introduced to adaptively weight and integrate all the features $A^{n,i}$. It should be noted that when a certain  data is constant or missing, the corresponding branch will be removed in DataInNet. For instance, in the Poisson equation, we will remove the branches corresponding to the initial conditions, namely $Q^{n,I}$ and $\hbar^{I}(\bar{\textbf{x}})$.

In certain PDEs, the magnitudes of specific physical quantities are typically extremely large. For instance, Poisson equations with high frequency solutions will result in source terms of considerable magnitude. Since we intend to incorporate such data into the network architecture, this may lead to problems such as gradient vanishing. Therefore, normalization of data is imperative. In this work, we use the following data processing method:
\begin{equation}\label{Nor}
\begin{aligned}
\bar{F}=\frac{F}{K^{*}|F|_{max}},
\end{aligned}
\end{equation}
where $K^{*}\geq 1$ is a manually adjustable hyperparameter. 

\begin{figure}[!htbp]
    \centering
    \includegraphics[height=8cm]{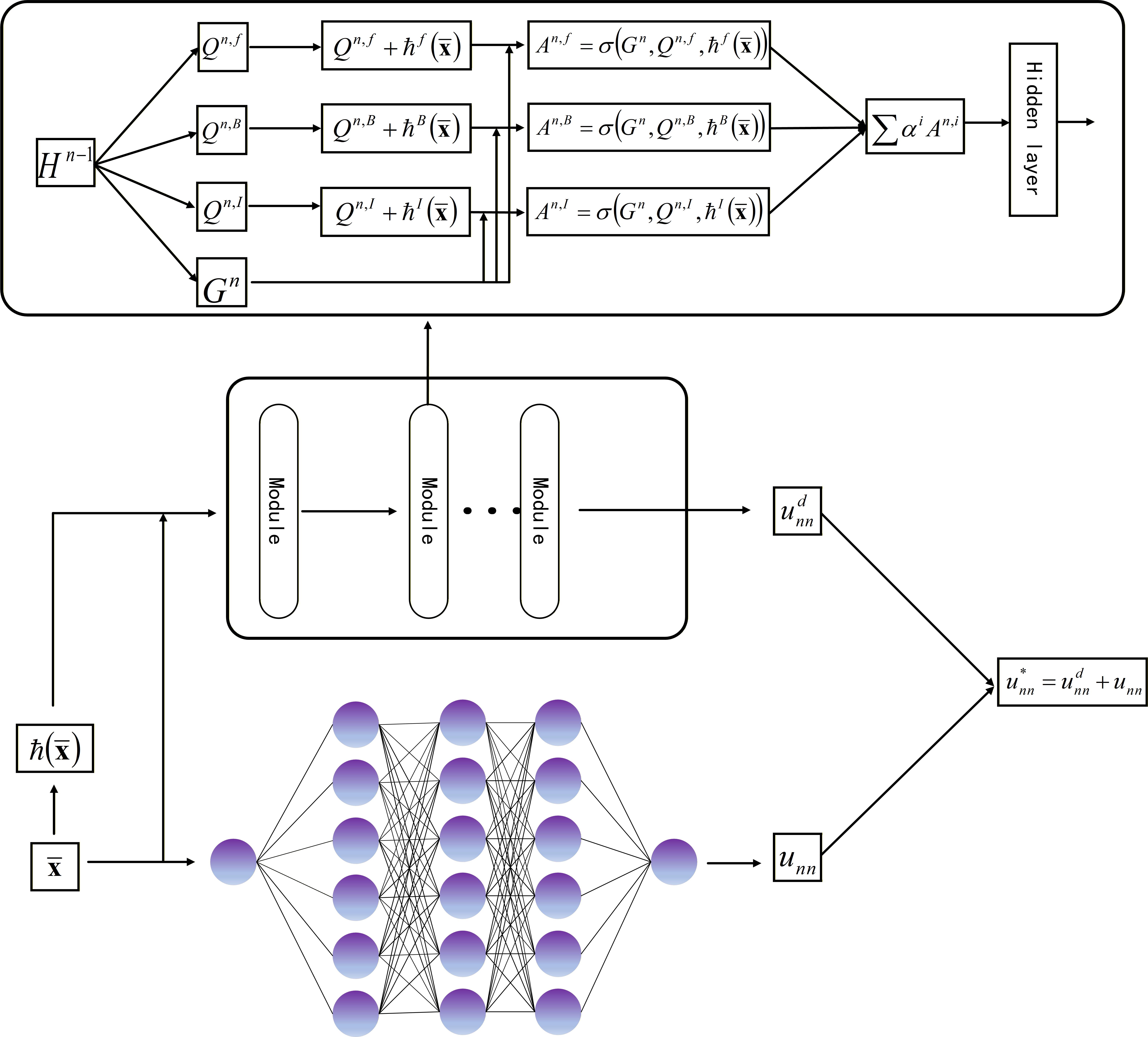}
    \caption{The structure of DataInNet.}
    \label{HLFNN}
    \end{figure}

We define the representation of DataInNet for problem (\ref{2.1}) as follows:
\begin{equation}\label{3.3}
\begin{aligned}
u_{nn}^{*}(\textbf{x},\theta)=u_{nn}(\textbf{x},t,\theta)+ u_{nn}^{d}(\textbf{x},t,\theta),
\end{aligned}
\end{equation}
where $u_{nn}(\cdot)$ and $u_{nn}^{d}(\cdot)$ represent the auxiliary neural network and data integration neural network, respectively. The network structure of DataInNet is shown in Figure \ref{HLFNN}. Then the loss function $L(\textbf{x},t,\theta)$ of problem (\ref{2.1}) is reformulated as:
\begin{equation}\label{3.4}
\begin{aligned}
L^{*}(\textbf{x},t,\theta)=&\tau\frac{1}{N}\sum\limits_{i=1}^{N} |D u_{nn}^{*}(\textbf{x}^{i},t^{i},\theta)- f(\textbf{x}^{i},t^{i})|^{2}\\ 
&+ \tau_{b}\frac{1}{N_{b}}\sum\limits_{i=1}^{N_{b}} |u_{nn}^{*}(\textbf{x}^{i}_{b},t^{i}_{b},\theta)-B(\textbf{x}^{i}_{b},t^{i}_{b})|^{2}\\
&+\tau_{0}\frac{1}{N_{0}}\sum\limits_{i=1}^{N_{0}}|u_{nn}^{*}(\textbf{x}^{i}_{0},0,\theta)-I(\textbf{x}^{i}_{0})|^{2}.\\
\end{aligned}
\end{equation}
By minimizing $L^{*}(\textbf{x},t,\theta)$, we obtain an approximation $u_{nn}^{*}(\textbf{x},t,\theta)$ of problem (\ref{2.1}). In order to obtain the optimal parameters $\theta^{*}=\{W^{*},b^{*}\}$, we use the gradient descent method to update the parameters. The update rules for step $k$ are as follows:
\begin{equation}\label{3.5}
\begin{aligned}
W^{k+1}=&W^{k}-\eta^{k} \frac{\partial L^{*}(\textbf{x},t,W^{k},b^{k})} {\partial W^{k}},\\
b^{k+1}=&b^{k}-\eta^{k} \frac{\partial L^{*}(\textbf{x},t,W^{k},b^{k})} {\partial b^{k}},\\
\end{aligned}
\end{equation}
where $\eta^{k}$ denotes the step size of the $k$-th iteration. 

\section{Numerical experiments}\label{section-3}
In this section, we investigate numerically the performance of DataInNet for solving PDEs. In all simulations, Adam optimizer with a learning rate of $0.001$ was used \cite{33}. The parameter $K^{*}$ is set to $2$. Furthermore,  The code accompanying this manuscript are publicly available at https://github.com/jczheng126/DataInNet. The relative $L^{2}$ error is used to evaluate the performance of the proposed method:
\begin{equation}\label{4.1}
\begin{aligned}
E_{L^2}= \frac{\sqrt{\sum\limits_{i=1}^{N}|u(\textbf{x}_{i})-u_{nn}^{*}(\textbf{x}_{i},\theta)|^{2}}} {\sqrt{\sum\limits_{i=1}^{N}|u(\textbf{x}_{i})|^{2}}},\\
\end{aligned}
\end{equation}
where $u(\textbf{x}_{i})$ denotes the true solution or the reference solution obtained using traditional numerical methods.

In DataInNet, a key design is to extend the scope of application of the boundary conditions and initial conditions from the geometric boundaries and initial time to the entire spatiotemporal computational domain. This means that for any input point within the domain, the model receives data on the boundary conditions and initial conditions calculated based on the coordinates of that point. This enables the model, when processing interior points, to explicitly perceive and incorporate the global physical constraints imposed by the boundary conditions and initial conditions. To systematically evaluate the effectiveness of this scope extension strategy, we consider a heat equation:
\begin{equation}\label{T1}
\left\{\begin{aligned}
u_{t}(x,t)=&\frac{1}{(10\pi)^{2}}u_{xx}(x,t),\qquad (x,t)\in [0,1]\times[0,1],\\
u(0,t)=&u(1,t)=0,\,\,\qquad\qquad t\in [0,1],\\
u(x,0)=&sin(10\pi x),\,\,\qquad\qquad x\in [0,1].\\
\end{aligned}
\right.
\end{equation}
The exact solution:
\begin{equation}\label{T2}
u(x,t) = e^{-t}sin(10\pi x).
\end{equation}

In this test, we design and compare two different strategies for integrating initial condition data:
\begin{itemize}
\item \textbf{Local-Input Strategy:} 
This strictly adheres to the classical PDE theory. The data of the initial condition is fed into model only when the temporal coordinate of an input point is exactly at the initial moment (i.e., $t = 0$.)

\item \textbf{Global-Input Strategy:}  
The initial condition is treated as a field defined over the entire spatial domain, and its scope of application is extended to the full spatiotemporal computational domain. For any spatiotemporal point within the computational domain (including interior and boundary points), the model receives the initial condition data corresponding to that point's spatial coordinates as part of its input.
\end{itemize}

We use $2,000$ interior points sampled by the Latin hypercube sampling (LHS) \cite{34}, $500$ initial points and $1,000$ boundary points as training points. A fully connected neural network with $4$ hidden layers serves as the auxiliary neural network. We use three modules with $48$ neurons to construct the data integration neural network. For each subnetwork, the $sin(x)$ function served as the activation function. The maximum number of iterations is set to $5,000$. 

\begin{figure}[!htbp]
    \centering
    \includegraphics[height=10cm]{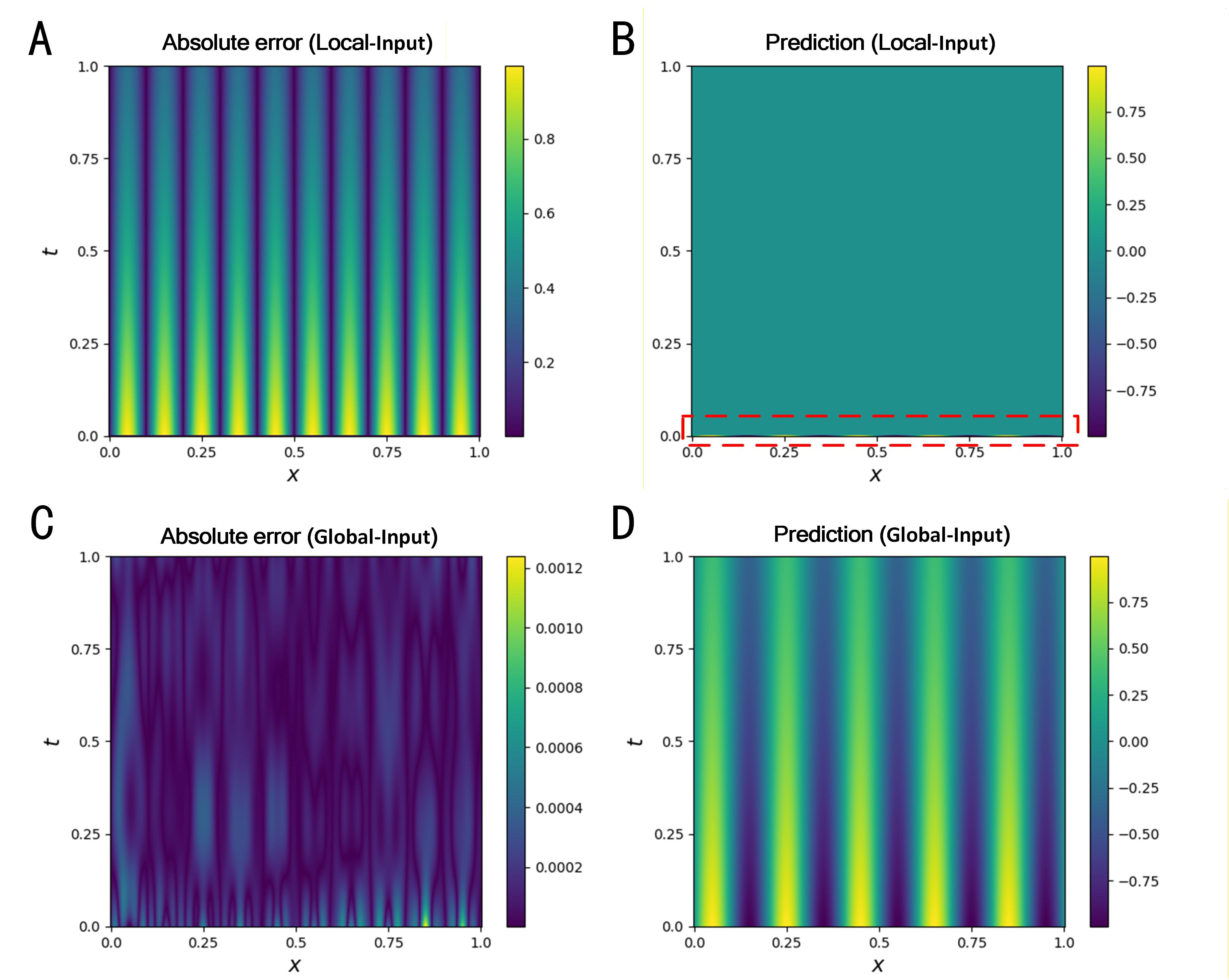}
    \caption{Heat equation (\ref{T1}).  \textbf{A}: Absolute error obtained using the local-input strategy. \textbf{B}: Prediction obtained using the local-input strategy. \textbf{C}: Absolute error obtained using the global-input strategy. \textbf{D}: Prediction obtained using the global-input strategy.}
    \label{GL}
    \end{figure}

Figure \ref{GL} (A) and (B) present the absolute error and the predicted solution obtained using the local-input strategy, respectively. The numerical results indicate that DataInNet based on this strategy fails to provide an accurate approximation. It is particularly noteworthy that, as observed in Figure \ref{GL} (B), the predicted solution is only accurate near the initial moment (t=0). This demonstrates that the local-input strategy cannot effectively propagate the physical constraints of the initial condition into the interior of the computational domain. In contrast, Figure \ref{GL} (C) and (D) display the results obtained using the global-input strategy. The model achieves significantly higher accuracy, with the maximum absolute error and relative $L^{2}$ error being only $1.22\times10^{-3}$ and $2.98\times10^{-4}$, respectively. Compared to the local-input strategy, the global-input strategy provides more effective guidance for DataInNet, making it the superior choice within this framework.

\begin{figure}[!htbp]
    \centering
    \includegraphics[height=8cm]{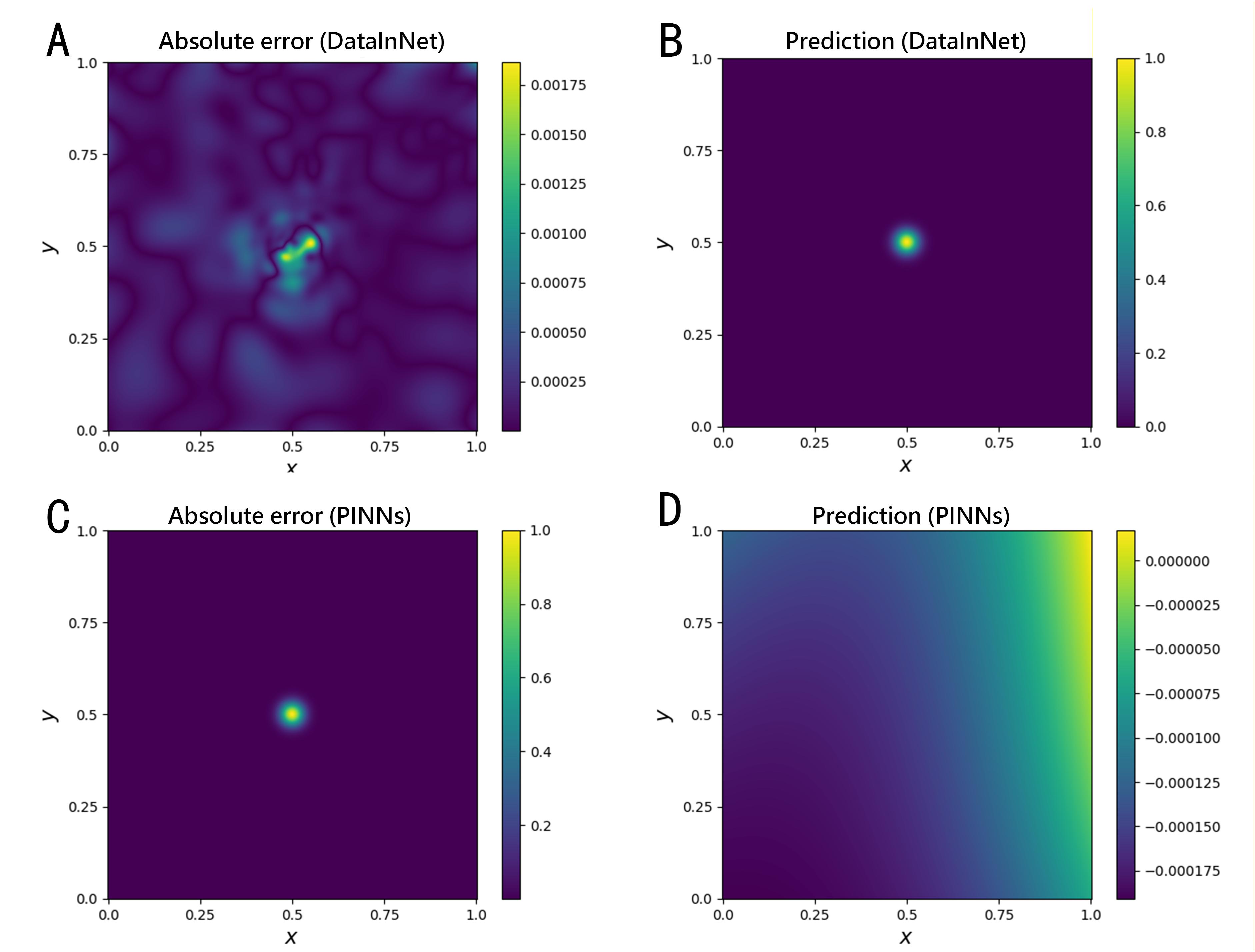}
    \caption{Example \ref{E1}: Poisson equation (\ref{e1}).  \textbf{A}: Absolute error for DataInNet. \textbf{B}: Prediction for DataInNet. \textbf{C}: Absolute error for PINNs. \textbf{D}: Prediction for PINNs.}
    \label{E1.1}
    \end{figure}

\begin{example}\label{E1}
We solve the two-dimensional Poisson equation using DataInNet:
\begin{equation}\label{e1}
\left\{\begin{aligned}
-\Delta u=&f(x,y),\qquad (x,y)\in \Omega,\\
u=&g(x,y),\qquad (x,y)\in \partial\Omega,\\
\end{aligned}
\right.
\end{equation}
where $\Omega=[0,1]\times[0,1]$, and the exact solution:
\begin{equation}\label{e1.1}
u(x,y) = e^{-1000((x-0.5)^2+(y-0.5)^2)}.
\end{equation}
\end{example}

We use $4,000$ interior points sampled by the LHS and $1,200$ boundary points as training points. A fully connected neural network with $4$ hidden layers serves as the auxiliary neural network. We use three modules with $48$ neurons to construct the data integration neural network. The weight $\tau_{b}=10,000$ and the maximum number of iterations is set to $70,000$. 

The solution (\ref{e1.1}) has a peak at $(0.5, 0.5)$ and decays rapidly away from $(0.5, 0.5)$. This may cause the neural network to converge to a trivial solution. Figure \ref{E1.1} (C) and (D) show the absolute error and the approximation obtained by PINNs, respectively. The numerical results indicate that PINNs fails to capture the sharp variation of the solution at $(0.5, 0.5)$. The absolute error obtained by DataInNet is shown in Figure \ref{E1.1} (A), where the maximum absolute error is $1.77 \times 10^{-3}$. Figure \ref{E1.1} (B) displays the predicted solution obtained by DataInNet. By integrating the source term and boundary conditions into the data integration neural network, DataInNet inherently excludes function classes that violate known constraints, thereby providing an accurate approximation. In Table \ref{E1:Tab1}, we record the performance of the model in various numbers of modules and training points. Additionally, we evaluated the influence of activation functions on numerical results. Table \ref{E1:Tab2} records the relative $L^{2}$ errors obtained by DataInNet with different activation functions. 

\begin{table}[htbp]
	\centering
	\caption{Relative $L^{2}$ errors in Example \ref{E1}}
	\label{E1:Tab1}
        \setlength{\tabcolsep}{6mm}
	\begin{tabular}{c|cccc ccc}
		\hline\hline\noalign{\smallskip}
		\diaghead{\theadfont Diag ColumnmnAB} {Interior points} {Modules} & 1 & 2 & 3 & 4 \\  
           \noalign{\smallskip}\hline\noalign{\smallskip}
            1000 & $9.27\times10^{-1}$ & $4.17\times10^{-1}$ & $3.17\times10^{-2}$ &  $1.79\times10^{-2}$ \\
            2000 & $7.64\times10^{-1}$ & $6.84\times10^{-2}$ & $8.12\times10^{-3}$  &  $9.23\times10^{-3}$  \\
            4000 & $3.44\times10^{-1}$ & $4.75\times10^{-2}$ & $4.07\times10^{-3}$ & $7.25\times10^{-3}$  \\
		\noalign{\smallskip}\hline
	\end{tabular}
\end{table}

\begin{table}[htbp]
	\centering
	\caption{Relative $L^{2}$ errors for different activation functions in Example \ref{E1}}
	\label{E1:Tab2}
        \setlength{\tabcolsep}{6mm}
        \renewcommand{\arraystretch}{1.5}
	\begin{tabular}{c|cccc ccc}
		\hline\hline\noalign{\smallskip}	
		{Activation functions} & sin(x) & cos(x) & tanh(x) & sigmoid(x)\\  
           \noalign{\smallskip}\hline\noalign{\smallskip}
            Relatvie $L^{2}$ errors & $4.03\times10^{-3}$ & $1.27\times10^{-2}$ & $1.43\times10^{-2}$ & $9.48\times10^{-2}$\\
		\noalign{\smallskip}\hline
	\end{tabular}
\end{table}

\begin{example}\label{E2}
We solve the two-dimensional Poisson equation defined on an L-shaped region:
\begin{equation}\label{e2}
\left\{\begin{aligned}
-\Delta u(x,y)=&f(x,y),\qquad (x,y)\in \Omega,\\
u(x,y)=&0, \,\,\qquad\qquad (x,y)\in \partial\Omega,\\
\end{aligned}\right.
\end{equation}
where $\Omega=(-1, 1)\times(-1,1)/[0,1)\times(-1,0]$. The exact solution is given as:
\begin{equation}
u(x,y)=\frac{3}{2}(1-x^{2})(1-y^{2})r^{\frac{2}{3}}sin(\frac{2\theta}{3}),
\end{equation}
where $r = \sqrt{x^2+y^2}$, $\theta=arctanh(\frac{y}{x})$.
\end{example}

In this test, $1,000$ interior points and $800$ boundary points are used for network training. We use a fully connected neural network with $4$ hidden layers as the auxiliary neural network, with $48$ neurons in each hidden layer. Four modules with $48$ neurons are used to construct the data integration neural network. For each subnetwork, the $tanh(x)$ function served as the activation function. The maximum number of iterations is set to $100,000$.

\begin{figure}[!htbp]
    \centering
    \includegraphics[height=8cm]{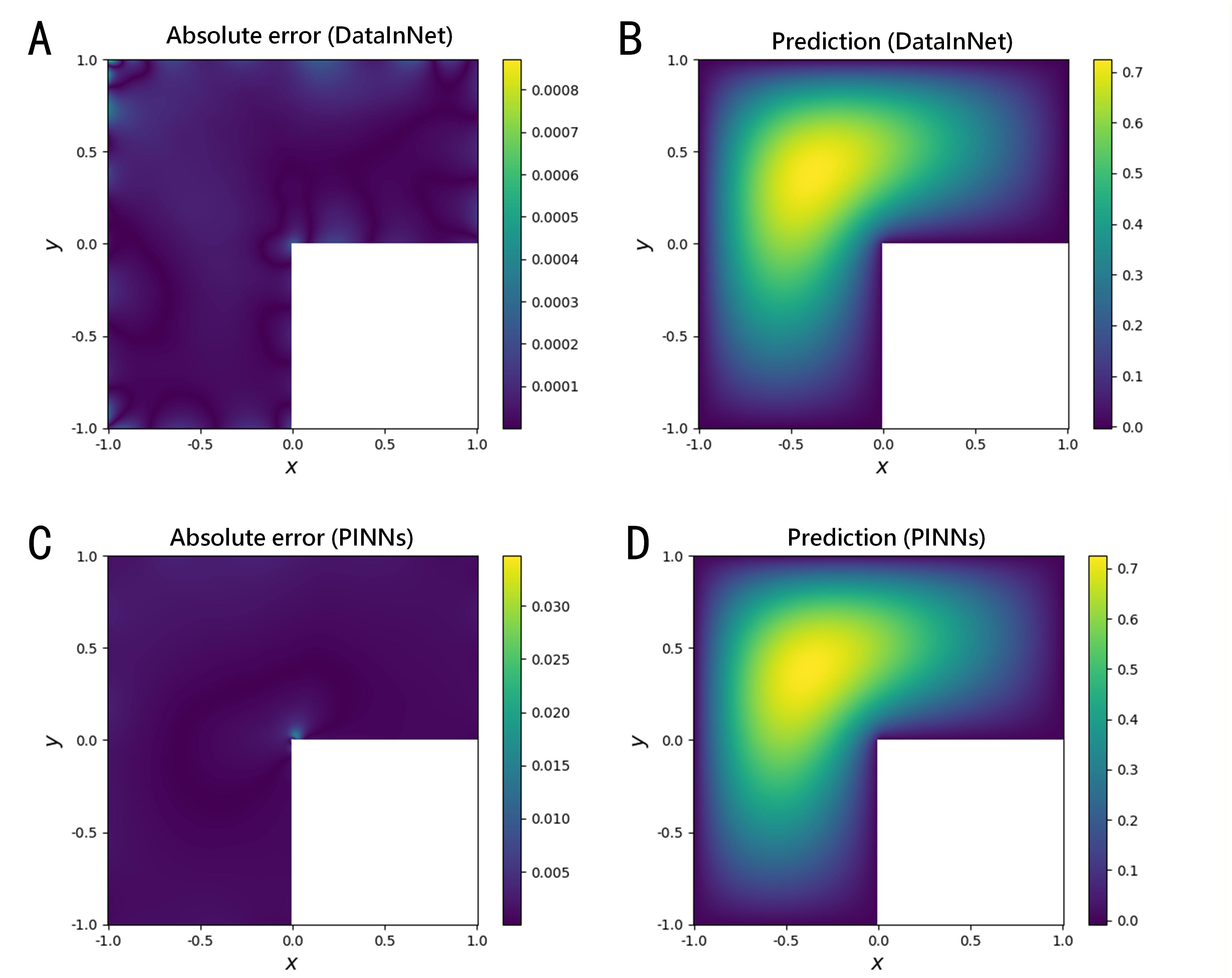}
    \caption{Example \ref{E2}: Poisson equation (\ref{e2}). \textbf{A}: Absolute error for DataInNet. \textbf{B}: Prediction for DataInNet. \textbf{C}: Absolute error for PINNs. \textbf{D}: Prediction for PINNs.}
    \label{E2.1}
    \end{figure}

Figure \ref{E2.1} (A) shows the absolute error for DataInNet. The maximum absolute error and relative $L^{2}$ error are $8.16\times10^{-4}$ and $1.71\times10^{-4}$, respectively. Figure \ref{E2.1} (B) shows the prediction for DataInNet. Figure \ref{E2.1} (C) and (D) display the absolute error and prediction of PINNs, respectively. The numerical results show that DataInNet provides an accurate approximation. It should be noted that, compared to PINNs, DataInNet achieves a significantly better approximation at the reentrant corner and yields lower errors.

 We test the performance of the model at different numbers of modules and training points. Table \ref{E2:Tab1} provides a detailed comparison. The numerical results in Table \ref{E2:Tab1} lead to a clear conclusion: As the number of modules used to build the data integration neural network increases and the number of sample points used for training grows, the accuracy of the approximated solution obtained by DataInNet continues to improve significantly. Table \ref{E2:Tab2} provides a detailed performance comparison of different activation functions under the same experimental configuration.

\begin{table}[htbp]
	\centering
	\caption{Relative $L^{2}$ errors in Example \ref{E2}}
	\label{E2:Tab1}
        \setlength{\tabcolsep}{6mm}
	\begin{tabular}{c|cccc ccc}
		\hline\hline\noalign{\smallskip}
		\diaghead{\theadfont Diag ColumnmnAB} {Interior points} {Modules} & 1 & 2 & 3 & 4 \\  
           \noalign{\smallskip}\hline\noalign{\smallskip}
            500 & $3.11\times10^{-1}$ & $7.32\times10^{-3}$ & $2.21\times10^{-3}$ &  $8.15\times10^{-4}$ \\
            1000 & $1.24\times10^{-1}$ & $2.31\times10^{-3}$ & $7.74\times10^{-4}$  &  $1.71\times10^{-4}$  \\
            2000 & $2.15\times10^{-2}$ & $1.45\times10^{-3}$ & $6.12\times10^{-4}$ & $1.43\times10^{-4}$  \\
		\noalign{\smallskip}\hline
	\end{tabular}
\end{table}

\begin{table}[htbp]
	\centering
	\caption{Relative $L^{2}$ errors for different activation functions in Example \ref{E2}}
	\label{E2:Tab2}
        \setlength{\tabcolsep}{6mm}
        \renewcommand{\arraystretch}{1.5}
	\begin{tabular}{c|cccc ccc}
		\hline\hline\noalign{\smallskip}	
		{Activation functions} & sin(x) & cos(x) & tanh(x) & sigmoid(x)\\  
           \noalign{\smallskip}\hline\noalign{\smallskip}
            Relatvie $L^{2}$ errors & $5.37\times10^{-4}$ & $6.59\times10^{-4}$ & $1.71\times10^{-4}$ & $2.51\times10^{-3}$\\
		\noalign{\smallskip}\hline
	\end{tabular}
\end{table}

\begin{example}\label{E3}
We consider the following Helmholtz equation:
\begin{equation}\label{e3}
\left\{\begin{aligned}
-\Delta u(x,y)-k^{2}u(x,y)=&f(x,y),\qquad (x,y)\in \Omega,\\
u(x,y)=&0, \,\,\qquad\qquad (x,y)\in \partial\Omega,\\
\end{aligned}\right.
\end{equation}
where $\Omega=[-1,1]\times[-1,1]$ and $k=1$. The exact solution is given as:
\begin{equation}
u(x,y)=sin(\pi x)sin(4\pi y).
\end{equation}
\end{example}

We employed LHS to generate $2,300$ training points, comprising $1,500$ interior points and $800$ boundary points. We utilized a fully connected neural network with $4$ hidden layers as the auxiliary neural network, each hidden layer containing $48$ neurons. Four modules, each with $48$ neurons, are used to construct the data integration neural network. For each subnetwork, the $sin(x)$ function served as the activation function, and the maximum number of iterations is set to $50,000$.

\begin{figure}[!htbp]
    \centering
    \includegraphics[height=8cm]{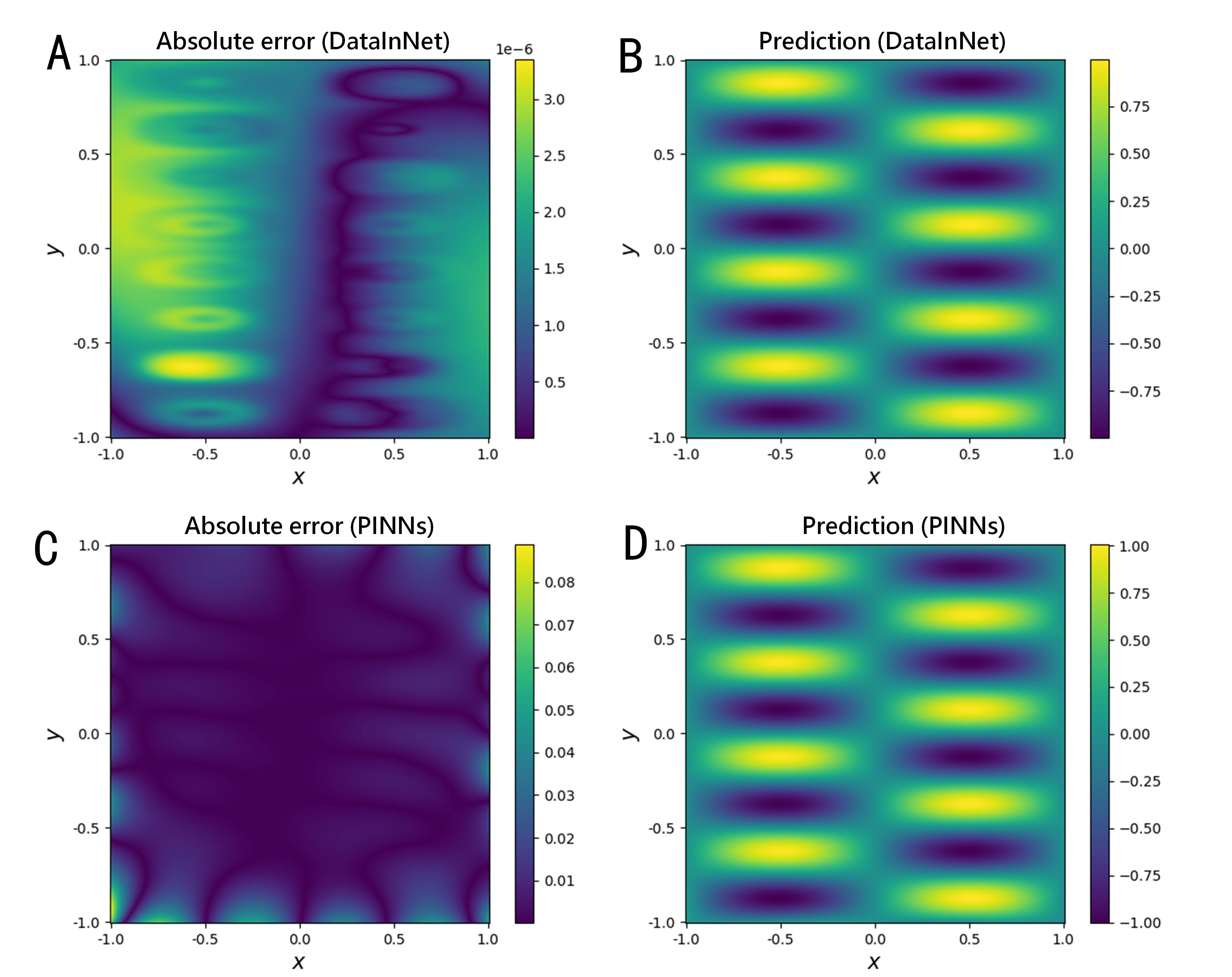}
    \caption{Example \ref{E3}: Helmholtz equation (\ref{e3}). \textbf{A}: Absolute error for DataInNet. \textbf{B}: Prediction for DataInNet. \textbf{C}: Absolute error for PINNs. \textbf{D}: Prediction for PINNs.}
    \label{E3.1}
    \end{figure}

Figure \ref{E3.1} (A) shows the absolute error for DataInNet. The maximum absolute error and relative $L^{2}$ error are $3.18\times10^{-6}$ and $3.23\times10^{-6}$, respectively. Figure \ref{E3.1} (B) shows the prediction for DataInNet. Figure \ref{E3.1} (C) and (D) display the absolute error and prediction of PINNs, respectively. To systematically evaluate the sensitivity of the DataInNet model performance to different activation functions, we conducted a series of controlled experiments. We tested the model's performance under varying numbers of modules and training points. Table \ref{E3:Tab1} shows detailed numerical results. As the number of modules used to build the data integration neural network increases and the number of sample points used for training grows, the accuracy of the approximation obtained by DataInNet is significantly improved. Table \ref{E3:Tab2} provides a detailed comparison. It can be observed that $sin(x)$ achieved the best approximation in this case. 
 
\begin{table}[htbp]
	\centering
	\caption{Relative $L^{2}$ errors in Example \ref{E3}}
	\label{E3:Tab1}
        \setlength{\tabcolsep}{6mm}
	\begin{tabular}{c|cccc ccc}
		\hline\hline\noalign{\smallskip}	
		\diaghead{\theadfont Diag ColumnmnAB} {Interior points} {Modules} & 1 & 2 & 3 & 4\\  
           \noalign{\smallskip}\hline\noalign{\smallskip}
            1000 & $1.14\times10^{-2}$ & $9.86\times10^{-3}$ &  $8.78\times10^{-5}$ & $8.28\times10^{-6}$\\
            1500 &  $2.03\times10^{-3}$  & $5.21\times10^{-4}$ & $2.51\times10^{-5}$ & $3.23\times10^{-6}$\\
            2000 & $5.97\times10^{-4}$ & $7.02\times10^{-5}$ & $3.75\times10^{-6}$  & $3.08\times10^{-6}$\\
		\noalign{\smallskip}\hline
	\end{tabular}
\end{table}

\begin{table}[htbp]
	\centering
	\caption{Relative $L^{2}$ errors for different activation functions in Example \ref{E3}}
	\label{E3:Tab2}
         \setlength{\tabcolsep}{6mm}
         \renewcommand{\arraystretch}{1.5}
	\begin{tabular}{c|cccc ccc}
		\hline\hline\noalign{\smallskip}	
		{Activation functions} & sin(x) & cos(x) & tanh(x) & sigmoid(x)\\  
           \noalign{\smallskip}\hline\noalign{\smallskip}
            Relatvie $L^{2}$ errors & $3.23\times10^{-6}$ & $8.69\times10^{-6}$ & $3.73\times10^{-5}$  & $7.18\times10^{-5}$\\
		\noalign{\smallskip}\hline
	\end{tabular}
\end{table}

\begin{figure}[!htbp]
    \centering
    \includegraphics[height=8cm]{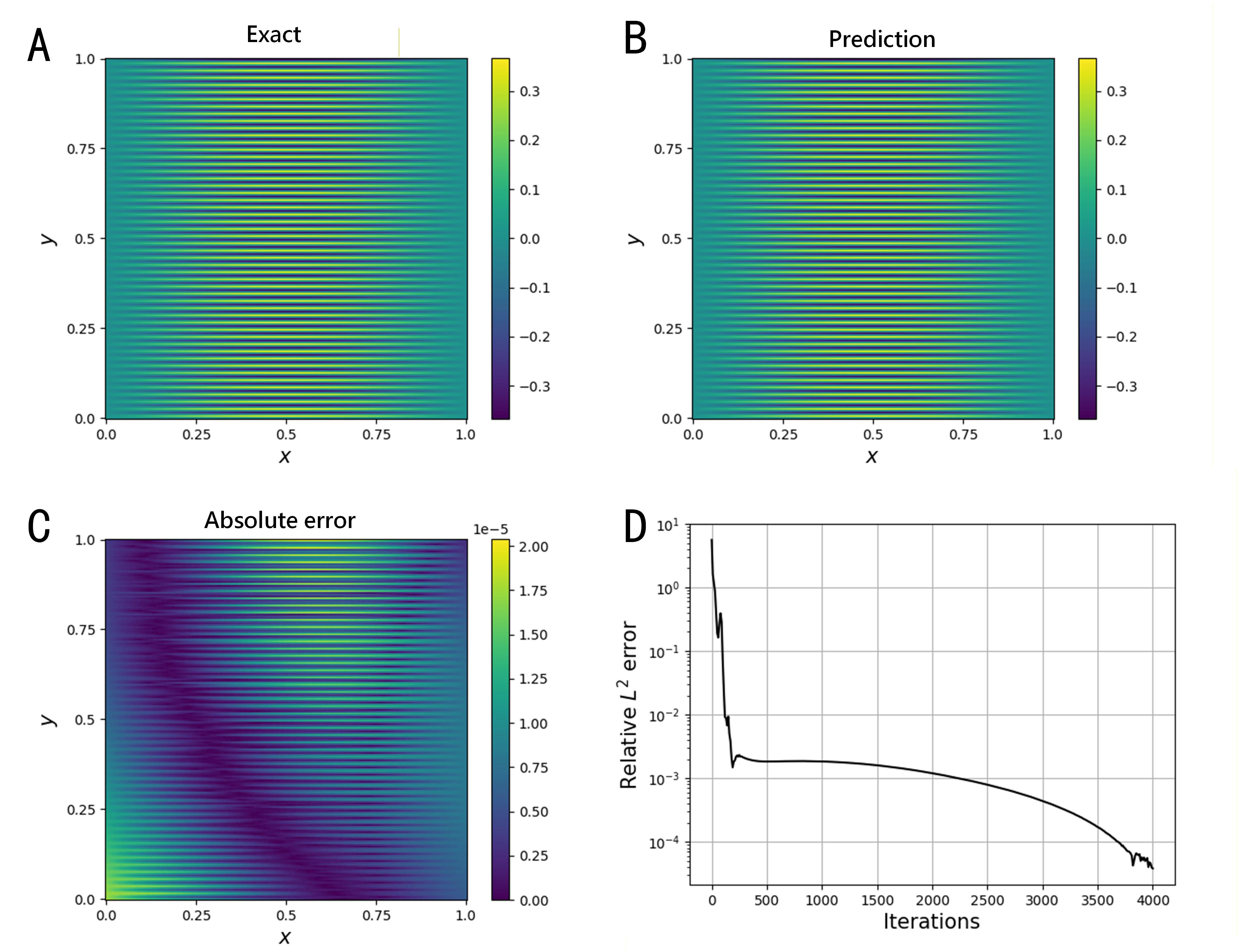}
    \caption{Example \ref{E4}: Wave equation (\ref{e4}). \textbf{A}: Exact solution (t=0.5). \textbf{B}: Predicte solution (t=0.5).
\textbf{C}: Absolute error. \textbf{D}: Relative $L^{2}$ error.}
    \label{E4.1}
    \end{figure}

\begin{table}[htbp]
	\centering
	\caption{Relative $L^{2}$ errors in Example \ref{E4}}
	\label{E4:Tab1}
        \setlength{\tabcolsep}{6mm}
	\begin{tabular}{c|cccc ccc}
		\hline\hline\noalign{\smallskip}	
		\diaghead{\theadfont Diag ColumnmnAB} {Interior points} {Modules} & 1 & 2 & 3 & 4\\  
           \noalign{\smallskip}\hline\noalign{\smallskip}
            1000 & $4.16\times10^{-2}$ & $1.77\times10^{-4}$ & $1.67\times10^{-4}$  & $7.64\times10^{-5}$\\
            2000 & $7.39\times10^{-3}$ & $1.17\times10^{-4}$ & $2.91\times10^{-4}$  & $3.81\times10^{-5}$\\
            3000 & $3.28\times10^{-3}$ & $1.72\times10^{-4}$ & $2.64\times10^{-4}$  & $1.71\times10^{-4}$\\
		\noalign{\smallskip}\hline
	\end{tabular}
\end{table}

\begin{table}[htbp]
	\centering
	\caption{Relative $L^{2}$ errors for different activation functions in Example \ref{E4}}
	\label{E4:Tab2}
        \setlength{\tabcolsep}{6mm}
        \renewcommand{\arraystretch}{1.5}
	\begin{tabular}{c|cccc ccc}
		\hline\hline\noalign{\smallskip}	
		{Activation functions} & sin(x) & cos(x) & tanh(x) & sigmoid(x)\\  
           \noalign{\smallskip}\hline\noalign{\smallskip}
            Relatvie $L^{2}$ errors & $1.41\times10^{-4}$ & $2.94\times10^{-4}$ & $1.15\times10^{-3}$ & $3.81\times10^{-5}$\\
		\noalign{\smallskip}\hline
	\end{tabular}
\end{table}

\begin{example}\label{E4}
We consider the following nonhomogeneous wave equation:
\begin{equation}\label{e4}
\left\{\begin{aligned}
u_{tt}(x,y,t)-\Delta u(x,y,t) + u_{t}(x,y,t)^{2}=&f(x,y,t),\qquad (x,y,t)\in \Omega \times [0, T],\\
u(x,y,t)=&g(x,y,t), \qquad (x,y)\in \partial\Omega, t\in [0,T],\\
u(x,y,0)=&h(x,y), \,\;\;\qquad(x,y)\in \Omega,\\
u_{t}(x,y,0)=&m(x,y),\,\;\qquad(x,y)\in \Omega,
\end{aligned}\right.
\end{equation}
where $\Omega=[0,1]^{2}$ and $T=1$, the exact solution:
\begin{equation}\label{e41}
u(x,y,t)=e^{-2t}sin(\pi x)sin(100\pi y).
\end{equation}
\end{example}

In this test, we use DataInNet to approximate solution (\ref{e41}) with high frequency components. Since the initial condition and source term provide the same frequency information, we only integrate the source term into the data integration neural network.

We employed LHS to generate $3,200$ training points, comprising $2,000$ interior points, $800$ boundary points and $400$ initial points. The weights $\tau=1$, $\tau_{b}=1,000$. We use a fully connected neural network with $4$ hidden layers as the auxiliary neural network, each layer containing $48$ neurons. Four modules, each with $48$ neurons, are used to construct the data integration neural network. The $sigmoid(x)$ function served as the activation function for each subnetwork and the maximum number of iterations is set to $4,000$.

Figure \ref{E4.1} (A) and (B) display the exact solution ($t=0.5$) and the predicted solution ($t=0.5$), respectively. It can be found that DataInNet provides an accurate approximation. The absolute error is presented in Figure \ref{E4.1} (C), where the maximum absolute error is $2.31\times10^{-5}$. The relative $L^{2}$ error is shown in Figure \ref{E4.1} (D). During the early stages of training, the relative $L^{2}$ error decreases rapidly to the order of $O(10^{-3})$. After 4000 iterations, DataInNet yields an approximation with relative $L^{2}$ error of $O(10^{-5})$. Furthermore, Table \ref{E4:Tab1} provides a detailed comparison concerning different numbers of modules and training points. As the number of modules and training points increases, the accuracy of the approximate solution obtained by the DataInNet is significantly improved. 
Table \ref{E4:Tab2} provides a detailed performance comparison of different activation functions under the same experimental configuration. The numerical results indicate that $sigmoid(x)$ is the most suitable and optimal activation function for this case. 

\begin{figure}[!h]
    \centering
    \includegraphics[height=8cm]{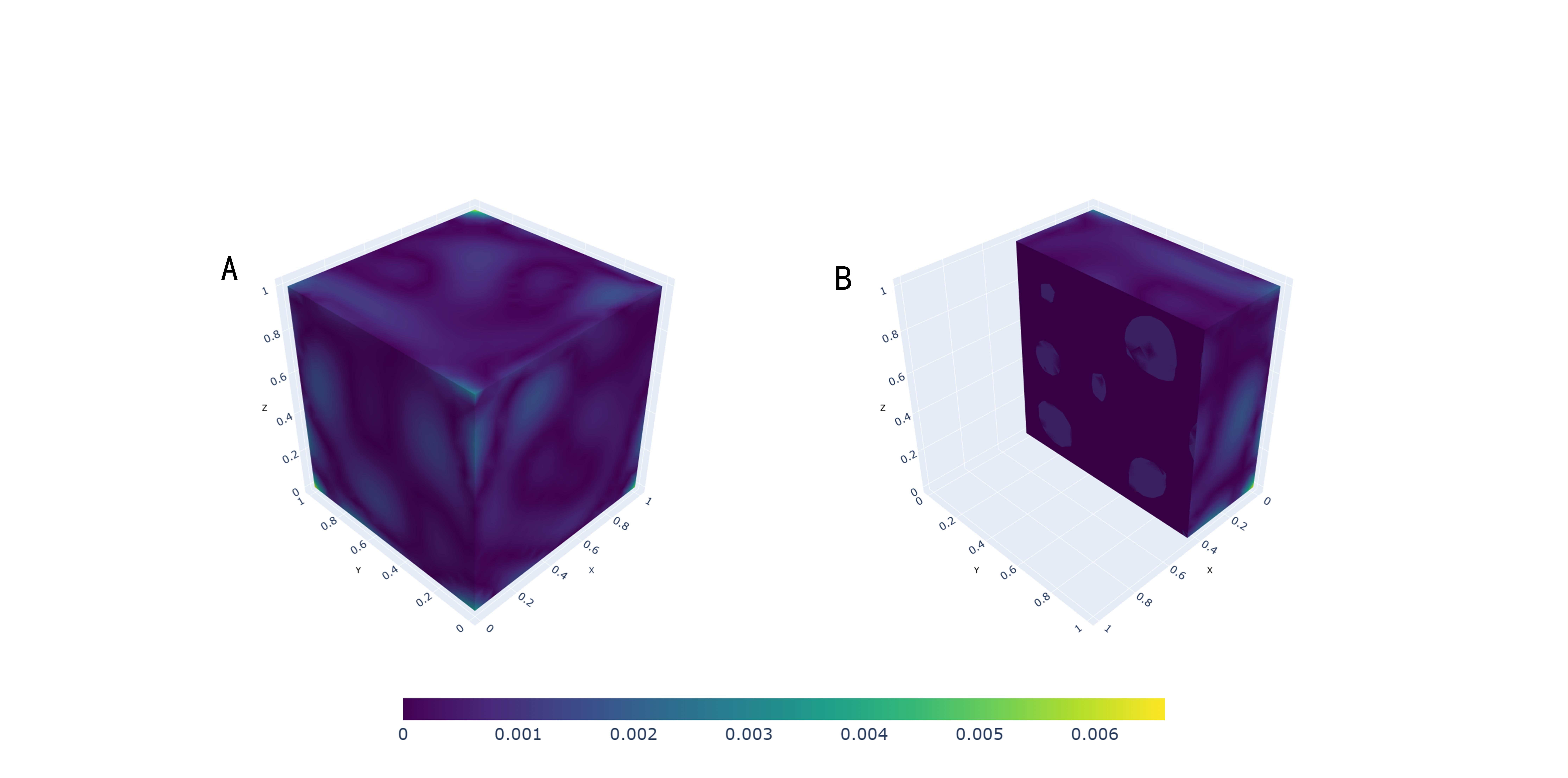}
    \caption{Example \ref{E5}: Three-dimensional Poisson equation (\ref{e5}). \textbf{A}: Absolute error at the boundary of the computational domain. \textbf{B}: Absolute error inside the computational domain.}
    \label{E5.1}
    \end{figure}

\begin{example}\label{E5}
Finally, we consider the following three-dimensional Poisson equation:
\begin{equation}\label{e5}
\left\{\begin{aligned}
-\Delta u(x,y,z)=&f(x,y,z),\qquad (x,y,z)\in \Omega,\\
u(x,y,z)=&g(x,y,z), \qquad(x,y,z)\in \partial\Omega,\\
\end{aligned}\right.
\end{equation}
where $\Omega=[0,1]^{3}$, and the exact solution:
\begin{equation}\label{e5.1}
\begin{split}
u(x,y,z)=sin(\pi x)sin(\pi y)sin(\pi z) + sin(10\pi x)cos(30\pi y)sin(20\pi z).
\end{split}
\end{equation}
\end{example}

In this test, we use LHS to generate $4,400$ training points, including $2,000$ interior points and $2,400$ boundary points. The data integration neural network consists of two modules with $48$ neurons, while the auxiliary neural network is a fully connected neural network with $4$ hidden layers. The maximum number of iterations is set to $40,000$.
 
\begin{figure}[!h]
    \centering
    \includegraphics[height=5cm]{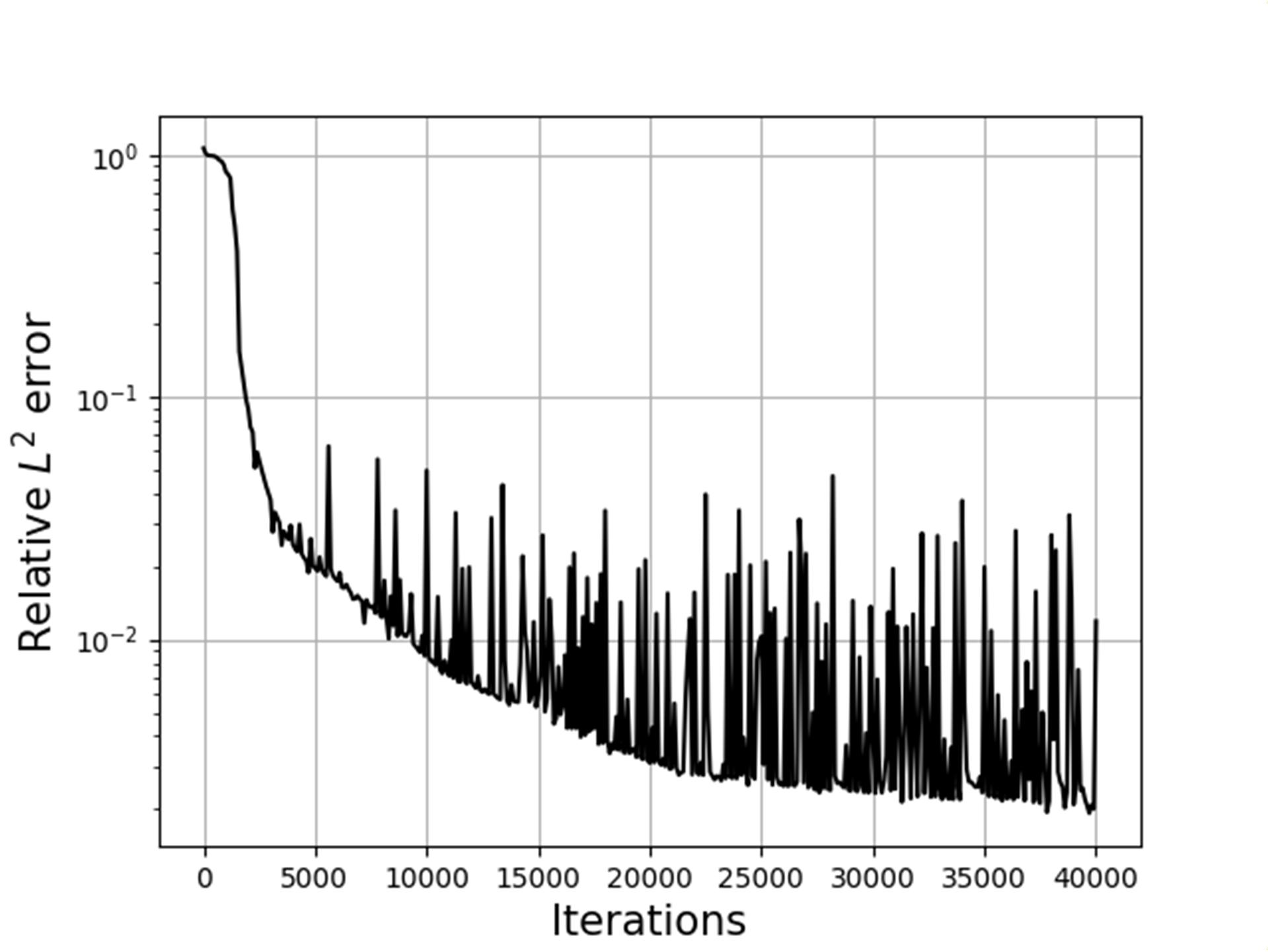}
    \caption{Example \ref{E5}: Three-dimensional Poisson equation (\ref{e5}). The decrease of relative $L^{2}$ error during the training process.}
    \label{E5.2}
    \end{figure}

Figure \ref{E5.1} displays the absolute error generated by DataInNet in solving problem (\ref{e5}), where the maximum absolute error and relative $L^{2}$ error are $6.08\times10^{-3}$ and $1.93\times10^{-3}$, respectively. Figure \ref{E5.2} illustrates the variation of the relative error during the training process. The numerical results demonstrate that DataInNet can provides an accurate approximation for three dimensional PDEs with high frequency components. Table \ref{E5:Tab1} provides a comparison concerning different numbers of modules and training points. Combined with previous examples, a series of results validates the effectiveness of DataInNet in solving PDEs. 
We also evaluated the performance of DataInNet in solving problem (\ref{e5}) using different activation functions. The numerical results are recorded in Table \ref{E5:Tab2}, where the best approximation is achieved using the $sin(x)$ activation function. 

\begin{table}[htbp]
	\centering
	\caption{Relative $L^{2}$ errors in Example \ref{E5}}
	\label{E5:Tab1}
        \setlength{\tabcolsep}{6mm}
	\begin{tabular}{c|cccc ccc}
		\hline\hline\noalign{\smallskip}	
		\diaghead{\theadfont Diag ColumnmnAB} {Interior points} {Modules} & 1 & 2 & 3 & 4\\ 
           \noalign{\smallskip}\hline\noalign{\smallskip}
            1000 & $7.95\times10^{-2}$ & $3.16\times10^{-3}$ & $2.93\times10^{-3}$  & $2.78\times10^{-3}$\\
            2000 & $9.48\times10^{-3}$ & $1.93\times10^{-3}$ & $2.56\times10^{-3}$  & $2.07\times10^{-3}$\\
            3000 & $8.71\times10^{-3}$ & $2.18\times10^{-3}$ & $2.36\times10^{-3}$  & $3.13\times10^{-3}$\\
		\noalign{\smallskip}\hline
	\end{tabular}
\end{table}

\begin{table}[htbp]
	\centering
	\caption{Relative $L^{2}$ errors for different activation functions in Example \ref{E5}}
	\label{E5:Tab2}
        \setlength{\tabcolsep}{6mm}
        \renewcommand{\arraystretch}{1.5}
	\begin{tabular}{c|cccc ccc}
		\hline\hline\noalign{\smallskip}	
		{Activation functions} & sin(x) & cos(x) & tanh(x) & sigmoid(x)\\  
           \noalign{\smallskip}\hline\noalign{\smallskip}
            Relatvie $L^{2}$ errors & $1.93\times10^{-3}$ & $4.86\times10^{-3}$ & $1.05\times10^{-2}$  & $8.67\times10^{-2}$\\
		\noalign{\smallskip}\hline
	\end{tabular}
\end{table}

\section{Summary and discussion}\label{section-4}
In this work, we propose data-integrated neural networks (DataInNet) for solving PDEs. DataInNet comprises two subnetworks: a data integration neural network dedicated to fusing various types of data, and an auxiliary neural network tasked with learning the residual physical information not captured by the data integration neural network. Notably, DataInNet inherently eliminates solutions that violate known  constraints, thereby guiding the optimization process toward more feasible function classes. This method reforms the conventional paradigm of physics-informed neural networks (PINNs)-where source terms, initial condition, and boundary conditions are merely imposed as soft constraints in the loss function-thus offering a novel perspective on the utilization of  prior information. It should be noted that when the source terms, initial conditions, and boundary conditions are all constant, there is no essential distinction between DataInNet and fully connected neural networks.

Numerical experiments demonstrate the effectiveness of DataInNet. These experiments encompass the Helmholtz equation, the Poisson equation defined on an L-shaped domain, and PDEs with high frequency solutions. Future research will focus on devising more efficient data integration methods and assessing the applicability of this approach in more challenging scenarios.


\section*{Acknowledgments}
This work was supported by NSFC Project (12431014), Project of Scientific Research Fund of the Hunan Provincial Science and Technology Department (2024ZL5017), 111 Project (No.D23017) and Guizhou Provincial Science and Technology Projects, China (No.QKHJC-ZK[2023]YB036).


\begin{thebibliography}{99}
\bibitem{1}
\textsc{L. Yang, D. Zhang, G.E. Karniadakis}, {Physics-informed generative adversarial networks for stochastic differential equations}, SIAM Journal on Scientific Computing, 42 (2020) A292-A317.
\bibitem{2}
\textsc{G. Pang, L. Lu, G.E. Karniadakis},  {fPINNs: fractional physics-informed neural networks}, SIAM Journal on Scientific Computing, 41 (2019) A2603-A2626.
\bibitem{3}
\textsc{H. Zhang, X. Shao, Z. Zhang, M. He},  {E-PINN: extended physics informed neural network for the forward and inverse problems of high-order nonlinear integro-differential equations}, International Journal of Computer Mathematics, 101 (7) (2024) 723-749.
\bibitem{4}
\textsc{Z. Wang, Z. Zhang}, {A mesh-free method for interface problems using the deep learning approach}, Journal of Computational Physics, 400 (2020) 108963.
\bibitem{5}
\textsc{S. Wang, P. Perdikaris}, {Deep learning of free boundary and Stefan problems}, Journal of Computational Physics, 428 (2021) 109914.
\bibitem{6}
\textsc{Y. Tseng, T. Lin, W. Hu, M. Lai},  {A cusp-capturing PINN for elliptic interface problems}, Journal of Computational Physics, 491 (2023) 112359.
\bibitem{7}
\textsc{A.K. Sarma, S. Roy, C. Annavarapu, P. Roy, S. Jagannathan},  {Interface PINNs (I-PINNs): A physics-informed neural networks framework for interface problems}, Computer Methods in Applied Mechanics and Engineering, 429 (2024) 117135.
\bibitem{8}
\textsc{X. Meng, Z. Li, D. Zhang, G.E. Karniadakis}, {PPINN: Parareal physics-formed neural network for time-dependent PDEs}, Computer Methods in Applied Mechanics and Engineering, 370 (2020) 113250.
\bibitem{9}
\textsc{A.D. Jagtap, E. Kharazmi, G.E. Karniadakis}, {Conservative physics-informed neural networks on discrete domains for conservation laws}, Computer Methods in Applied Mechanics and Engineering, 365 (2020) 113028.
\bibitem{10}
\textsc{Y. Tseng, T. Lin, W. Hu, M. Lai},  {A cusp-capturing PINN for elliptic interface problems}, Journal of Computational Physics, 491 (2023) 112359.
\bibitem{11}
\textsc{E. Kharazmi, Z. Zhang, G.E. Karniadakis}, {hp-VPINNs: Variational physics-informed neural networks with domain decomposition}, Computer Methods in Applied Mechanics and Engineering, 374 (2021) 113547.
\bibitem{12}
\textsc{Y. Yang, P. Perdikaris}, {Adversarial uncertainty quantification in physics-informed neural networks}, Journal of Computational Physics, 394 (2019) 136-152.
\bibitem{13}
\textsc{A.D. Jagtap, G.E. Karniadaki}, {Extended physics-informed neural networks (XPINNs): A generalized space-time domain decomposition based deep learning framework for nonlinear partial differential equations}, Communications in Computational Physics, 28 (2020) 2002-2041.
\bibitem{14}
\textsc{X. Jin, S. Cai, H. Li, G.E. Karniadakis},  {NSFnets (Navier-Stokes flow nets): Physics-informed neural networks for the incompressible Navier-Stokes equations}, Journal of Computational Physics, 426 (2021) 109951.
\bibitem{15}
\textsc{W. Hu, T. Lin, M Lai},  {A discontinuity capturing shallow neural network for elliptic interface problems}, Journal of Computational Physics, 469 (2022) 111576.
\bibitem{16}
\textsc{S. Wang, S. Sankaran, P. Perdikaris},  {Respecting causality for training physics-informed neural networks},  Computer Methods in Applied Mechanics and Engineering, 412 (2024) 116813.
\bibitem{17}
\textsc{Z. Xiang, W. Peng, X. Liu, W. Yao},  {Self-adaptive loss balanced physics-informed neural networks}, Neurocomputing, 496 (2022) 11-34.
\bibitem{18}
\textsc{E. Kharazmi, Z. Zhang, G.E. Karniadakis}, {VPINNs: Variational physics-informed neural networks for solving partial differential equations}, Preprint at arXiv. abs/1912.00873 (2019).
\bibitem{19}
\textsc{Z.-Q.J. Xu, Y. Zhang, T. Luo, Y. Xiao, Z. Ma},  {Frequency principle: Fourier analysis sheds light on deep neural networks}, Communications in Computational Physics, 28 (5) (2020) 1746-1767.
\bibitem{20}
\textsc{Z.-Q.J. Xu, L. Zhang, W. Cai},  {On understanding and overcoming spectral biases of deep neural network learning methods for solving PDEs}, Journal of Computational Physics, 530 (2025) 113905.
\bibitem{21}
\textsc{L. Lu, X. Meng, Z. Mao, and G.E. Karniadakis},  {DeepXDE: A deep learning library for solving differential equations}, SIAM Review, 63 (1) (2021) 208-228.
\bibitem{22}
\textsc{W. E, C. Ma, L. Wu},  { Machine learning from a continuous viewpoint, I}, Science China Mathematics, 63 (2020) 2233-2266.
\bibitem{23}
\textsc{M. Tancik, P. Srinivasan, B. Mildenhall, S. Fridovich-Keil, N. Raghavan, U. Singhal, R. Ramamoorthi, J. Barron, R. Ng}, {Fourier features let networks learn high frequency functions in low dimensional domains}, Advances in Neural Information Processing Systems, 33 (2020) 7537–7547.
\bibitem{24}
\textsc{S. Wang, X. Yu, P. Perdikaris}, {When and why PINNs fail to train: A neural tangent kernel perspective}, Journal of Computational Physics, 449 (2022) 110768.
\bibitem{25}
\textsc{S. Wang, H. Wang, P. Perdikaris}, {On the eigenvector bias of Fourier feature networks: From regression to solving multi-scale PDEs with physics-informed neural networks}, Computer Methods in Applied Mechanics and Engineering, 384 (2021) 113938.
\bibitem{26}
\textsc{L. Lu, X. Meng, Z. Mao, G.E. Karniadakis},  {DeepXDE: A deep learning library for solving differential equations}, SIAM Review, 63 (1) (2021) 208-228.
\bibitem{27}
\textsc{H. Fan, Z. Tan}, {Novel and general discontinuity-removing PINNs for elliptic interface problems}, Journal of Computational Physics, 529 (2025) 113861.
\bibitem{28}
\textsc{Y. Wang, C. Yao}, {Multi-stage neural networks: Function approximator of machine precision}, Journal of Computational Physics, 504 (2024) 112865.
\bibitem{29}
\textsc{J. Huang, R. You, T. Zhou}, {Frequency-adaptive multi-scale deep neural networks}, Computer Methods in Applied Mechanics and Engineering, 437 (2025) 117751. 
\bibitem{30}
\textsc{A.D. Jagtap, K. Kawaguchi, G.E. Karniadakis}, {Adaptive activation functions accelerate convergence in deep and physics-informed neural networks}, Journal of Computational Physics, 404 (2020) 109136.
\bibitem{31}
\textsc{F. Cao, X. Guo, X. Dong, D. Yuan}, {wbPINN: Weight balanced physics-informed neural networks for multi-objective learning}, Applied Soft Computing, 170 (2025) 112632.
\bibitem{32}
\textsc{S. Wang, Y. Teng, P. Perdikaris}, {Understanding and mitigating gradient pathologies in physics-informed neural networks}, SIAM Journal on Scientific Computing, 43 (5) (2021) A3055–A3081.
\bibitem{33}
\textsc{D.P. Kingma, J. Ba},  {Adam: A method for stochastic optimization}, Preprint at arXiv: 1412.6980 (2014).
\bibitem{34}
\textsc{M. Stein},  {Large sample properties of simulations using Latin hypercube sampling}, Technometrics, 29 (2) (1987) 143-151.

\end{thebibliography}
\end{document}